\documentclass[12pt]{article}
\usepackage{amsmath,amsfonts,amssymb}
\begin{document}
\newtheorem{teo}{Theorem}
\newtheorem{cor}{Corollary}
\newtheorem{lm}{Lemma}
\newtheorem{rem}{Remark}
\newtheorem{pr}{Proposition}
\renewcommand{\theteo}{\arabic{pr}.}
\renewcommand{\theteo}{\arabic{teo}.}
\renewcommand{\thecor}{\arabic{cor}.}
\renewcommand{\thelm}{\arabic{lm}.}
\renewcommand{\therem}{\arabic{rem}.}
\newcommand{\indlim}{\operatornamewithlimits{ind\ lim}}
\renewcommand{\refname}{\begin{center} \mdseries {\sf \large References}
\end{center}}
\begin{center}{\bf ON EXISTENCE OF BOUNDARY VALUES OF POLYHARMONIC FUNCTIONS}
%\end{center}
\footnote{
Supported by CRDF and Ukr. Government (Project UM 1-2507-OD-03)}
\end{center}
\vspace*{3mm}
\begin{center} {\textmd \sf M.L. GORBACHUK and S.M. TORBA}
\end{center}
\vspace*{3mm}
\begin{quote} ABSTRACT. \ In trigonometric series terms all polyharmonic
functions inside the unit disk are described. For such functions it is proved
the existence of their boundary values on the unit circle in the space of
hyperfunctions. The necessary and sufficient conditions are presented for
the boundary value to belong to certain subspaces of the space of
hyperfunctions.
\end{quote}
\vspace*{3mm}

The purpose of this paper is to find necessary and sufficient conditions for
a solution of the equation ${\Delta}^mu = 0$ inside a domain to have a
limit on the boundary of the domain in various functional spaces. We consider
the simplest situation where a domain is the unit disk $K = \{z = re^{it}, \ 0
\le r < 1, \ 0 \le t \le 2\pi\}$. The case of $m = 1$ has been investigated
during 20th century by a lot of mathematicians (we refer for details to
[1 - 5]. The case of $m = 2$ was considered in [6]. For an arbitrary $m$ the
problem of existence of boundary values in the space $L_2(\partial K) \
(\partial K$ - is the unit circle) was discussed in [7].

{\bf 1.} \ Denote by $D = D(\partial K)$ the set of all infinitely
differentiable functions on $\partial K$. We say that a sequence
${\varphi}_n \in D$ converges to $\varphi \in D, \ n \to \infty$, and write
${\varphi}_n \stackrel{D}{\to}\varphi$, if for every $k \in {\Bbb N}_0 =
\{0, 1, 2, \dots\}$, the sequence ${\varphi}_n^{(k)}(t)$ converges
to ${\varphi}^{(k)}(t)$ uniformly in $t \in \partial K$. Let also ${\frak A =
\frak A}(\partial K)$ be the set of analytic functions on $\partial K$.
The convergence in $\frak A$ is introduced in the following way: a sequence
${\varphi}_n \in \frak A$ converges to $\varphi$ in $\frak A \ ({\varphi}_n
\stackrel{\frak A}{\to}\varphi)$ if there exists a neighbourhood $U$ of
$\partial K$ in which all the functions ${\varphi}_n(t)$ converge to
$\varphi(t)$ uniformly on any compact set from $U$.

For a number $\alpha > 0$ we put
$$
{\frak A}_{\alpha} = \{\varphi \in D \bigl| \exists c > 0 \ \forall {\Bbb N}_0 \
\max\limits_{t \in \partial K}|{\varphi}^{(k)}(t)| \le c{\alpha}^kk!\}.
$$
The linear set ${\frak A}_{\alpha}$ is a Banach space with respect to the norm
$$
\|\varphi\|_{{\frak A}_{\alpha}} = \sup\limits_{k \in {\Bbb N}_0}
\frac{\max\limits_{t \in \partial K}|{\varphi}^{(k)}(t)|}{{\alpha}^kk!}.
$$
It is not hard to show that if $\alpha < \alpha'$, then ${\frak A}_{\alpha}
\subseteq {\frak A}_{\alpha'}$,
$$
{\frak A} = \indlim\limits_{\alpha \to \infty}{\frak A}_{\alpha},
$$
and the dense continuous embeddings
$$
{\frak A} \subset D \subset L_p(\partial K) = L_p, \quad 1 \le p < \infty,
$$
hold.

Let $D'$ and ${\frak A}'$ are the spaces of continuous antilinear
functionals on $D$ (distributions) and $\frak A$ (hyperfunctions),
respectively (see [8]). In the following,
$<F, \varphi>$ denotes an action of the functional $F \in {\frak A}' \
(F \in D')$ onto $\varphi \in {\frak A} \ (\varphi \in D)$. By convergence
in ${\frak A}'$ (in $D'$) we mean the weak one, that is, $F_n
\stackrel{{\frak A}'}{\to}F \ (F_n \stackrel{D'}{\to}F)$ if for any
$\varphi \in {\frak A} \ (\varphi \in D)$, the number sequence
$<F_n, \varphi>$ converges to $<F, \varphi>$.

As $e_k(t) = e^{ikt} \in {\frak A} \ (k \in \Bbb Z)$, the Fourier coefficients
$c_k(F) = <F, e_k>$ can be determined for $F \in {\frak A}'$. It is known
(see e.g. [9]) that
$$
\sum\limits_{k = -n}^{n} c_k(F)e^{ikt} \stackrel{{\frak A}'}{\to} F,
$$
and one can easily verify that the below assertion is valid.
\begin{pr} \
The following equivalence relations hold:
$$
\begin{array}{lcl}
F \in D & \Longleftrightarrow & \forall \alpha > 0 \
\exists c > 0 \quad |c_k(F)| \le c|k|^{-\alpha}; \\
F \in {\frak A} & \Longleftrightarrow & \exists \alpha > 0 \
\exists c > 0 \quad |c_k(F)| \le ce^{-\alpha |k|}; \\
F \in {D'} & \Longleftrightarrow & \exists \alpha > 0 \
\exists c > 0 \quad |c_k(F)| \le c|k|^{\alpha}; \\
F \in {\frak A}' & \Longleftrightarrow & \forall \alpha > 0 \
\exists c > 0 \quad |c_k(F)| \le ce^{\alpha |k|}.
\end{array}
$$
Moreover, the series $\sum\limits_{k = -\infty}^{\infty} c_k(F)e^{ikt}$
converges to $F$ in the corresponding space. The sequence $\{F_n\}_{n \in \Bbb N}$, whose
elements $F_n$ belong to one of the spaces $D, {\frak A}, D'$ or ${\frak A}'$,
converges to $F$ in this space if and only if the constants $c$ and $\alpha$ in
the above estimates for $|c_k(F_n)|$ do not depend on $n$ and for any $k \in
{\Bbb Z}, \ c_k(F_n) \to c_k(F), \ n \to \infty$.
\end{pr}

{\bf 2.} \ A function $u(r, t) = u(re^{it}) \in C^{2m}(K)$ is called
$m$-harmonic in $K$ if it satisfies the equation
$$
{\Delta}^m u(r,t) = 0, \quad 0 \le r < 1, \ t \in [0, 2\pi]. \eqno (1)
$$
Note, that no conditions on the behaviour of $u(r, t)$ near $\partial K$ are
imposed.
\begin{teo} \
In order that a function $u(r, t) \in C^{2m}(K)$ be $m$-harmonic in $K$, it
is necessary and sufficient that the representation
$$
u(r,t) = \sum\limits_{j = 1}^{m}(r^2 - 1)^{j - 1}
\sum\limits_{k = -\infty}^{\infty} c_k(F_j)r^{|k|}e^{ikt}, \quad F_j \in
{\frak A}', \eqno (2)
$$
be admissible, where $F_j$ are uniquely determined by $u(r, t)$.
\end{teo}
{\it Proof}. \ By Proposition 1,
$$
\forall \alpha > 0 \
\exists c_j = c_j(\alpha) > 0 \ \forall k \in {\Bbb Z} \quad
|c_k(F_j)| \le c_je^{\alpha |k|}.
$$
So the series $\sum\limits_{k = -\infty}^{\infty} c_k(F_j)r^{|k|}e^{ikt}$
converges uniformly in the disk $\overline{K_R} =
\{z \in {\Bbb C}: |z| \le R\}$ of radius $R < e^{-\alpha}$ and determines an
infinitely differentiable function there. The direct check shows that
the functions
$(r^2 - 1)^{j-1}\sum\limits_{k = -\infty}^{\infty} c_k(F_j)r^{|k|}e^{ikt}, \ j =
1, 2, \dots, m,$ \ satisfy (1) in $K_R$. Since $\alpha > 0$ is arbitrary,
these functions are solutions of the equation (1) inside $K$.

To prove the necessity, suppose at first $m = 1$. Let $u(r, t)$ be a harmonic
function in $K$. Then for a fixed $r < 1, \ u(r, t)$ is infinitely
differentiable in $t$, and it may be written in the form
$$
u(r, t) = \sum\limits_{k = -\infty}^{\infty} c_k(r)e^{ikt}, \quad
c_k(r) = \frac{1}{2\pi}\int\limits_{0}^{2\pi} u(r, t)e^{-ikt}\,dt, \eqno (3)
$$
where the series and all its derivatives converge uniformly in
$t \in [0, 2\pi]$. The coefficients $c_k(r)$ are infinitely differentiable on
[0, 1) and satisfy the equation
$$
r^2{c_k}''(r) + r{c_k}'(r) = k^2c_k(r).
$$
Hence,
$$
c_k(r) = r^{|k|}c_k, \quad c_k \in \Bbb C.
$$
It follows from the convergence of the series in (3) that
$$
\forall r < 1 \quad  r^{|k|}|c_k| = e^{-\alpha |k|}|c_k| \le c,
$$
where $\alpha = -\ln r > 0$ is arbitrary. By Proposition 1, $c_k$ are the
Fourier coefficients of a certain hyperfunction $F_1$, and
$$
u(r, t) = \sum\limits_{k = -\infty}^{\infty} r^{|k|}c_k(F_1)e^{ikt}, \quad
F_1 \in {\frak A}'. \eqno (4)
$$
Thus, the representation (2) is valid when $m = 1$.

Assume the representation (2) to be true for an $(m - 1)$-harmonic
inside $K$ function $u(r, t) \ (m \ge 2)$, and we shall prove that such
a representation holds for an $m$-harmonic function.

If $u(r, t)$ is an $m$-harmonic function, then $\Delta u (r, t)$ is an
$(m - 1)$-harmonic one. By assumption, there exist $E_j \in {\frak A}', \ j =
1, 2, \dots, m - 1,$ such that
$$
\Delta u (r, t) = \sum\limits_{j = 1}^{m - 1}(r^2 - 1)^{j - 1}
\sum\limits_{k = -\infty}^{\infty} c_k(E_j)r^{|k|}e^{ikt}. \eqno (5)
$$
If we choose $\widetilde{u} \in C^2(K)$ so that
$$
\Delta (u (r, t) - \widetilde{u}(r, t)) = 0, \eqno (6)
$$
then, because of (4), we shall have
$$
u(r, t) = \widetilde{u}(r, t) +
\sum\limits_{k = -\infty}^{\infty}r^{|k|}c_k(F_1)e^{ikt}, \quad
F_1 \in {\frak A}'.
$$

Let us find at first $\widetilde{u}(r, t)$ in the case where the equation (5)
is of the form
$$
\Delta u (r, t) = \sum\limits_{k = -\infty}^{\infty}r^{|k|}c_k(E_1)e^{ikt},
\quad E_1 \in {\frak A}'.
$$
By using the identity
$$
\Delta\left((r^2 - 1)^j \sum\limits_{k = -\infty}^{\infty}r^{|k|}c_k(F)e^{ikt}\right) =
4j \sum\limits_{k = -\infty}^{\infty}\big[(r^2 - 1)^{j - 1}(|k| + j) +
(j - 1)(r^2 - 1)^{j - 2}\big]r^{|k|}c_k(F)e^{ikt}
$$
for $F \in {\frak A}'$, one can verify that the function ${\widetilde{u}}_2(r, t) = u_2(r, t)$, where
$$
u_2(r, t) = \frac{1}{4}(r^2 - 1)\sum\limits_{k = -\infty}^{\infty}
\frac{r^{|k|}}{|k| + 1}c_k(E_1)e^{ikt}
$$
$({\widetilde{u}}_1(r, t) \equiv 0)$, satisfies (6). Set $c_k(F_2) =
\frac{c_k(E_1)}{4(|k| + 1)}$. By Proposition 1,
$F_2 \in {\frak A}'$. So, in the case under consideration
$$
u(r, t) = \sum\limits_{j = 1}^2(r^2 - 1)^{j - 1}
\sum\limits_{k = -\infty}^{\infty} r^{|k|}c_k(F_j)e^{ikt}.
$$

Suppose now that we know solutions $u_l(r, t)$ of the equations
$$
\Delta u (r, t) = \sum\limits_{j = 1}^l(r^2 - 1)^{j - 1}
\sum\limits_{k = -\infty}^{\infty} r^{|k|}c_k(E_j)e^{ikt}, \eqno (8)
$$
for all $l \le s, \ s \le m - 2$ is fixed. We show how to find a solution
of the equation
$$
\Delta u (r, t) = \sum\limits_{j = 1}^{s + 1}(r^2 - 1)^{j - 1}
\sum\limits_{k = -\infty}^{\infty} r^{|k|}c_k(E_j)e^{ikt}. \eqno (9)
$$
We put
$$
u_{s + 2}(r, t) = (r^2 - 1)^{s + 1}
\sum\limits_{k = -\infty}^{\infty}
\frac{r^{|k|}}{4(s + 1)(|k| + s + 1)}c_k(E_{s + 1})e^{ikt}.
$$
It follows from (7) and (8) that if $u(r, t)$ is a solution of (9), then
$$
\Delta(u - u_{s + 2})(r, t) = \sum\limits_{j = 1}^{s + 1}(r^2 - 1)^{j - 1}
\sum\limits_{k = -\infty}^{\infty}r^{|k|}c_k(E_j)e^{ikt} -
$$
$$
\sum\limits_{k=-\infty}^{\infty}(r^2-1)^sr^{|k|}c_k(E_{s+1})e^{ikt}-\sum\limits_{k = -\infty}^{\infty}\left[\frac{s(r^2 - 1)^{s - 1}}{|k| + s + 1}
+ (r^2 - 1)^s\right] r^{|k|}c_k(E_{s + 1})e^{ikt} =
$$
$$
\sum\limits_{j = 1}^{s - 1}(r^2 - 1)^{j - 1}
\sum\limits_{k = -\infty}^{\infty} r^{|k|}c_k(E_j)e^{ikt} +
(r^2 - 1)^{s - 1}\sum\limits_{k = -\infty}^{\infty}r^{|k|}\left[c_k(E_s) -
\frac{s\cdot c_k(E_{s + 1})}{|k| + s + 1}\right]e^{ikt}.
$$
Taking into account that $E_s, E_{s + 1} \in {\frak A}'$, we conclude, by
Proposition 1, that there exists $E'_s \in {\frak A}'$ such that
$$
c_k(E'_s) = c_k(E_s) - \frac{s\cdot c_k(E_{s + 1})}{|k| + s + 1},
$$
whence
$$
\Delta(u - u_{s + 2})(r, t) =  \sum\limits_{j = 1}^{s}(r^2 - 1)^{j - 1}
\sum\limits_{k = -\infty}^{\infty} r^{|k|}c_k(E'_j)e^{ikt}, \quad E'_j
\in {\frak A}',
$$
where $E'_j = E_j$ as $j = 1, \dots, s - 1$. By assumption, we can find
${\widetilde{u}}_{s +1}(r, t)$ so that
$$
\Delta(u - u_{s + 2} - {\widetilde{u}}_{s +1})(r, t) = 0.
$$
Setting
$$
{\widetilde{u}}_{s +2}(r, t) = u_{s + 2}(r, t) +
{\widetilde{u}}_{s +1}(r, t),
$$
we arrive at the equality
$$
\Delta(u - {\widetilde{u}}_{s +2})(r, t) = 0.
$$

It is not hard to observe that for the desired function
$\widetilde{u}(r, t)$ we have the formula
$$
\widetilde{u}(r, t) = {\widetilde{u}}_m(r, t) = u_m(r, t) + u_{m- 1}(r, t)+
\dots + u_2(r, t) =
$$
$$
(r^2 - 1)^{m - 1}\sum\limits_{k = -\infty}^{\infty}
\frac{r^{|k|}}{4(m - 1)(|k| + m - 1)}c_k(E_{m - 1})e^{ikt} + \dots
$$
$$
+ (r^2 - 1)^2 \sum\limits_{k = -\infty}^{\infty}
\frac{r^{|k|}}{4\cdot 2(|k| + 2)}c_k(E_2)e^{ikt} + (r^2 - 1)
\sum\limits_{k = -\infty}^{\infty}
\frac{r^{|k|}}{4(|k| + 1)}c_k(E_1)e^{ikt}.
$$
Then
$$
u(r, t) =  \sum\limits_{j = 1}^{m}(r^2 - 1)^{j - 1}
\sum\limits_{k = -\infty}^{\infty} r^{|k|}c_k(F_j)e^{ikt}, \quad F_j
\in {\frak A}',
$$
where
$$
c_k(F_j) = \frac{c_k(E_{j - 1})}{4(j - 1)(|k| + j - 1)}.
$$

Since for $F \in {\frak A}'$
$$
r^{|k|}c_k(F) \to c_k(F), \ r \to 1, \ \mbox{\rm and} \quad |r^{|k|}c_k(F)| <
|c_k(F)|,
$$
we have, by Proposition 1, that
$$
(r^2 - 1)^{j - 1}\sum\limits_{k = -\infty}^{\infty} r^{|k|}c_k(F_j)e^{ikt}
\stackrel{{\frak A}'}{\to}
\left\{ \begin{array}{lcl}
F_1 & \mbox{\rm if} & j = 1 \\
0 & \mbox{\rm if} & j > 1,
\end{array}
\right.
$$
as $r \to 1$. The elements $F_j \in {\frak A}'$ are determined uniquely
by the function $u(r, t)$ in the following way:
$$
F_1 = \lim\limits_{r \to 1} u(r, \cdot), \quad F_{j + 1} =
\lim\limits_{r \to 1} \frac{u(r, \cdot) -  \sum\limits_{p = 1}^{j}(r^2 - 1)^{p - 1}
\sum\limits_{k = -\infty}^{\infty} r^{|k|}c_k(F_p)e^{ik\cdot}}{(r^2 - 1)^j},
$$
where the limit is taken in the space ${\frak A}'$. This completes the proof.

Because of harmonicity in $K$ of the functions
$$
u_j (r, t) = \sum\limits_{k = -\infty}^{\infty} r^{|k|}c_k(F_j)e^{ikt},
$$
the representation (2) implies, in particular, the next assertion (cf. [4]).
\begin{cor} \
Let $u(r, t)$ be an $m$-harmonic in $K$ function. Then it admits a representation
of the form
$$
u(r, t) =  \sum\limits_{j = 1}^{m}(r^2 - 1)^{j - 1}u_j(r, t), \eqno (10)
$$
where the functions $u_j(r, t)$ are harmonic in $K$.
\end{cor}
When proving the theorem, it was also established the following fact.
\begin{cor} \
If $u(r, t)$ is an $m$-harmonic in $K$ function, then there exists its radial
boundary value $u(1, \cdot)$ on $\partial K$ in the space ${\frak A}'$, that is,
$$
u(r, \cdot) \stackrel{{\frak A}'}{\to} u(1, \cdot) \quad \mbox{\rm as} \ r \to 1.
$$
\end{cor}

{\bf 3.} \  Let $\Phi$ be a complete linear Hausdorff space such that
the continuous embeddings
$$
{\frak A} \subset \Phi \subset {\frak A}'
$$
hold. We say that $F \in \Phi$ is a boundary value on $\partial K$ of an
$m$-harmonic in $K$ function $u(r, t)$ and write $F = u(1, \cdot)$ if
$u(r, \cdot) \stackrel{\Phi}{\to} F$ as $r \to 1$.

It is seen from Theorem 1 and Corollary 2 that every $m$-harmonic in $K$
function has a boundary value in ${\frak A}'$. Moreover, each element
$F \in {\frak A}'$  is the boundary value of a certain $m$-harmonic in $K$
function. The natural question arises: under what conditions on an
$m$-harmonic in $K$ function $u(r, t)$ its boundary value $u(1, \cdot)$
belongs to $\Phi$?
\begin{teo} \
The  boundary value $u(1, \cdot)$ of an $m$-harmonic in $K$ function $u(r, t)$
belongs to the space $\Phi$ if and only if the set $\{u(r, \cdot)\}_{r < 1}$
is compact in $\Phi$.
\end{teo}
{\it Proof}. {\it Necessity.} \ It is known that if $r_0 < 1$, then
$u(r_0, \cdot) \in \frak A$, and $u(r, \cdot) \stackrel{\frak A}{\to}
u(r_0, \cdot)$ as $1 >r \to r_0$. Since the embedding ${\frak A} \subset \Phi$
is continuous, $u(r, \cdot) \stackrel{\Phi}{\to} u(r_0, \cdot) \ (r \to r_0)$.
By assumption, $u(r, \cdot) \stackrel{\Phi}{\to}u(1, \cdot)$ if $r \to 1$.
So, the set $\{u(r, \cdot)\}_{r < 1}$ is compact in $\Phi$.

{\it Sufficiency.} \ Let the set $\{u(r, \cdot)\}_{r < 1}$ be compact in
$\Phi$. Suppose $r \to 1$. Then there exists a subsequence $r_k \to 1$ such
that $u(r_k, \cdot)$ converges in $\Phi \ (r_k \to 1)$ to a certain
element $F \in \Phi$. Since $\Phi \subset {\frak A}$ continuously,
$u(r_k, \cdot)$ converges in ${\frak A}'$. Taking into account that
$u(r, \cdot) \stackrel{{\frak A}'}{\to} u(1, \cdot)$ as $r \to 1$, we have
$u(1, \cdot) = F \in \Phi$ which completes the proof.

In the partial case where $\Phi = L_2(\partial K)$, Theorem 2 was obtained in
[7]. By using compactness criteria for sets, one can find the sufficient
conditions for the boundary value of a polyharmonic function to belong to
$L_p(\partial K), \ 1 \le p < \infty$. For instance, the following assertion
is valid.
\begin{cor} \
Let $u(r, t)$ be an $m$-harmonic inside the disk $K$ function. In order that $u(r, t)$
have a boundary value in $L_p = L_p(\partial K)$, it is necessary and sufficient
that:

1) \ $\sup\limits_{0 \le r < 1}\|u(r, \cdot)\|_{L_p} < \infty$;

2) \ $\int\limits_0^{2\pi} \big|u\big(re^{i(t - \tau)}\big) - u(re^{it})\big|^p\,dt \to 0
\ (\tau \to 0) \mbox{\rm\ uniformly in} \ r \in [0, 1)$.
\end{cor}

Now we consider in more detail the case of $L_2$. Let
$$
{\frak B}_j = \Bigg\{F \in {\frak A}'\Biggl| \|F\|_{{\frak B}_j} =
\sup\limits_{0 \le r < 1} (1 - r^2)^j\left(\sum\limits_{k = -\infty}^{\infty}
r^{2|k|}|c_k(F)|^2\right)^{1/2} < \infty\Bigg\}.
$$
The set ${\frak B}_j$ with norm $\|\cdot\|_{{\frak B}_j}$ forms a Banach space.
\begin{teo} \
If $u(r, t)$ is an $m$-harmonic in $K$ function, then
$$
\sup\limits_{0 \le r < 1} \int\limits_0^{2\pi} |u(r, t)|^2\,dt < \infty
\Longleftrightarrow F_1 \in L_2, \ F_j \in {\frak B}_{j - 1} \
\mbox{\rm if} \ 2 \le j \le m,
$$
where $F_j$ are taken from representation (2). Moreover,
$u(r, \cdot) \to F_1 \ (r \to 1)$ weakly in the space $L_2$.
\end{teo}
{\it Proof.} \ Assume that in the representation (2) $F_1 \in L_2,
\ F_j \in {\frak B}_{j - 1} \ (j =2, \dots, m)$. Then
$$
\int\limits_0^{2\pi} |u(r, t)|^2\,dt = \sum\limits_{k = -\infty}^{\infty}
r^{2|k|}\big|c_k(F_1) + (r^2 - 1)c_k(F_2) + \dots + (r^2 - 1)^{m - 1}c_k(F_m)\big|^2 \le
$$
$$
m \sum\limits_{j = 1}^{m} \sum\limits_{k = -\infty}^{\infty}
r^{2|k|}(r^2 - 1)^{2(j - 1)}|c_k(F_j)|^2 \le c.
$$

Conversely, let $\sup\limits_{0 \le r < 1} \|u(r, \cdot)\|_{L_2} < \infty$.
Then, as was shown in [4, Lemma 7], each summand in (10) is bounded, too:
$$
\sup\limits_{0 \le r < 1} \|(1 - r^2)^{j - 1}u_j(r, \cdot)\|_{L_2} < \infty,
\quad j = 1, 2, \dots, m.
$$
This is equivalent to the inequality
$$
\sup\limits_{0 \le r < 1} (1 - r^2)^{2(j - 1)}
\sum\limits_{k = -\infty}^{\infty} r^{2|k|}|c_k(F_j)|^2 < \infty,
$$
that is, $F_1 \in L_2,
\ F_j \in {\frak B}_{j - 1} \ (j =2, \dots, m)$.

It still remains to prove the weak convergence of $u(r, \cdot)$ to
$F_1 \ (r \to 1)$ in $L_2$. Since $u(r, \cdot) \stackrel{{\frak A}'}{\to}F_1$
as $r \to 1$ and $e^{ikt} \in {\frak A} \ (k \in \Bbb Z)$, we have
$$
\lim\limits_{r \to 1} \int\limits_0^{2\pi} u(r, t) e^{ikt}\,dt =
\lim\limits_{r \to 1} r^{|k|}[c_{-k}(F_1) + (r^2 - 1)c_{-k}(F_2) + \dots +
c_{-k}(F_m)] =
$$
$$
c_k(F_1) = <F_1, e_k> =
\int\limits_0^{2\pi} F_1(t) e^{ikt}\,dt.
$$
Thus, $u(r, \cdot) \to F_1 \ (r \to 1)$ weakly in $L_2$ on a total set, and
$\sup\limits_{0 \le r < 1}\|u(r, \cdot)\|_{L_2} < \infty$. It follows from here
that $u(r, \cdot) \to F_1 \ (r \to 1)$ weakly in $L_2$. The proof is
complete.

Let $u(r, t)$ be a harmonic in $K$ function. It follows from (2) that
$$
\|u(r, \cdot)\|_{L_2}^2 = \sum\limits_{k = -\infty}^{\infty} r^{2|k|}
|c_k(F_1)|^2.
$$
In view of $\|u(r, \cdot)\|_{L_2} \le c$, the well-known Fatou lemma and
the Lebesgue theorem on passage to the limit yield
$$
\sum\limits_{k = -\infty}^{\infty}|c_k(F_1)|^2 < \infty, \quad F_1 = u(1, \cdot)
\in L_2, \quad \|u(r, \cdot)\|_{L_2} \to \|u(1, \cdot)\|_{L_2}, \ r \to 1.
$$
Therefore the weak convergence of $u(r, t)$ to $u(1, t)$ implies the strong one.
As was shown in [7], in the case of $m = 2$ the boundedness of
$\|u(r, \cdot)\|_{L_2}$ does not guarantee the convergence of $u(r, \cdot) \
(r \to 1)$ in $L_2$.

We pass now to the Sobolev spaces
$$
W_2^{\alpha} =  W_2^{\alpha}(\partial K) = \big\{F \in {\frak A}'\big|
\sum\limits_{k = -\infty}^{\infty}|k|^{2\alpha}|c_k(F)|^2 < \infty\big\}, \
\alpha \in \Bbb R.
$$
The following statement is valid.
\begin{teo} \
The embeddings
$$
W_2^{-j} \subset {\frak B}_j \subset W_2^{-j - 0} =
\bigcap\limits_{\varepsilon > 0} W_2^{-j - \varepsilon}
$$
hold.
\end{teo}
{\it Proof}. \ Since the function $f(r) = (1 - r^2)^{2j}r^{2k}, \ j,k \in
{\Bbb N}_0,$ reaches its maximum at the point $r^2 = \frac{k}{k + 2j}$, and
$$
\max\limits_{0 \le r < 1}f(r) = \left(\frac{2j}{k + 2j}\right)^{2j}
\left(\frac{k}{k + 2j}\right)^{k} < \frac{c_j}{k^{2j}},
$$
we have
$$
\sup\limits_{0 \le r < 1} (1 - r^2)^{2j}\sum\limits_{k = -\infty}^{\infty}
r^{2|k|}|c_k(F)| <
c_j \sum\limits_{k = -\infty}^{\infty}\frac{|c_k(F)|^2}{|k|^{2j}},
$$
that is, $W_2^{-j} \subset {\frak B}_j$.

Suppose now $F \in {\frak B}_j$. Then, substituting $z:=r^2$,
$$
\exists c > 0 \quad (1 - z)^{2j}\sum\limits_{k = -\infty}^{\infty}
z^{|k|}|c_k(F)|^2 < c.
$$
Multiplying this inequality by $(1 - z)^{-\alpha}$ and then integrating
along $[0, 1)$, we obtain
$$
\sum\limits_{k = -\infty}^{\infty} |c_k(F)|^2
\int\limits_{0}^{1} (1 - z)^{2j - \alpha}z^{|k|}\,dz < \infty. \eqno (11)
$$
If we put $\delta = 2j -\alpha + 1$, we get for $n \in \Bbb N$
$$
a_n = \int\limits_{0}^{1} (1 - z)^{2j - \alpha}z^{n}\,dz =
\frac{(n)!}{\delta(\delta + 1)\dots(\delta + n)}.
$$
Since
$$
\frac{a_n}{a_{n + 1}} = 1 + \frac{\delta}{n} + O\left(\frac{1}{n^2}\right),
$$
the relation $a_n = O\left(\frac{1}{n^{\delta}}\right)$ is fulfilled. Taking
in (11) $\alpha = 1 - \varepsilon, \ \varepsilon \in (0, 1)$, we conclude
that
$$
\sum\limits_{k = -\infty}^{\infty}\frac{|c_k(F)|^2}{|k|^{2j + \varepsilon}}
< \infty,
$$
that is, $F \in W_2^{j - \varepsilon}$, which completes the proof.

The next theorem is devoted to the question on the existence of boundary
values in the space $D'$ of distributions.
\begin{teo} \
In order that an $m$-harmonic in $K$ function $u(r, t)$ admit a representation
of the form (2) with $F_j \in D' \ (j = 1, \dots, m)$, it is necessary and
sufficient that
$$
\exists \alpha \ge 0 \ \exists c > 0 \ \
\sup\limits_{t \in [0, 2\pi]} |u(re^{it})| \le c(1 - r)^{-\alpha}.
\eqno (12)
$$
\end{teo}
{\it Proof.} \ Let the inequality (12) hold. Then for $p \in {\Bbb N}, \ p >
\alpha$, the function $v(r, t) = (1 - r^2)^pu(r, t)$ is $(m + p)$-harmonic
in $K$, and it is not difficult to verify that
$$
\sup\limits_{0 \le r < 1} \|v(r, \cdot)\|_{L_2} \le 2\pi c.
$$
By Theorems 3,4, the function $v(r, t)$ may be represented in the form (2)
where $F_j \in W_2^{-(j + p)}$. Since $D' =
\bigcup\limits_{\alpha > 0}W_2^{-\alpha}$, we have $F_j \in D'$.

The necessity of condition (12) for $m = 1$ was proved in [5]. Namely,
it was shown there that for a harmonic function of the form
$$
u(r, t) = \sum\limits_{k = -\infty}^{\infty}r^{|k|}c_k(F)e^{ikt}, \quad f \in
D'
$$
there exists $\alpha \ge 0$ such that
$$
|u(r, t)| \le c(1 - r)^{-\alpha}.
$$
If we take $\alpha = \max\limits_{j} {\alpha}_j$, where ${\alpha}_j$ corresponds
to $F_j$ from (2), we obtain the estimate (12) for an $m$-harmonic
function ($m$ is arbitrary).
\begin{cor} \
An $m$-harmonic in $K$ function $u(r, t)$ has a boundary value in $D'$ if
and only if it satisfies (12).
\end{cor}

For a number $\beta > 1$  we put
$$
{\frak G}_{\{\beta\}} = {\frak G}_{\{\beta\}}(\partial K) =
\{\varphi \in D\bigl| \exists \alpha > 0 \ \exists c > 0 \
\forall k \in {\Bbb N}_0 \quad \max\limits_{t \in \partial K}|{\varphi}^{(k)}(t)| \le
c{\alpha}^kk^{k\beta}\}. \eqno (13)
$$
The linear space ${\frak G}_{\{\beta\}}$  is endowed with the inductive limit
topology of the Banach spaces ${\frak G}_{\{\beta, \alpha\}}$ of functions
$\varphi \in D$  satisfying (13) with a fixed constant $\alpha$. The norm in
${\frak G}_{\{\beta, \alpha\}}$ is defined as
$$
\|\varphi\|_{{\frak G}_{\{\beta, \alpha\}}} =
\sup\limits_{k \in {\Bbb N}_0} \frac{\max\limits_{t \in \partial K}
|{\varphi}^{(k)}(t)|}{{\alpha}^kk^{k\beta}}.
$$
It is evident, that
$$
{\frak A} \subset {\frak G}_{\{\beta\}} \subset D \subset L_2 \subset D'
\subset {\frak G}'_{\{\beta\}} \subset {\frak A}',
$$
where ${\frak G}'_{\{\beta\}}$ denotes the dual of ${\frak G}_{\{\beta\}}$.
\begin{teo} \
An $m$-harmonic in $K$ function $u(r, t)$ admits a representation of
the form (2) with $F_j \in {\frak G}'_{\{\beta\}} \ (j = 1, \dots, m)$ if and
only if
$$
\forall \alpha > 0 \ \exists c = c(\alpha) > 0 \quad
\sup\limits_{0 \le r < 1}|u(r, t)| \le ce^{\alpha(1 - r)^{-q}}, \quad
q = \frac{1}{\beta - 1}.       \eqno (14)
$$
\end{teo}

The proof follows the scheme like that in Theorem 5 if to take into account
that
$$
F \in {\frak G}'_{\{\beta\}} \Longleftrightarrow \forall \alpha > 0 \quad
|c_k(F)| < ce^{-\alpha|k|^{1/\beta}},
$$
and the series $\sum\limits_{k = -\infty}^{\infty} c_k(F)e^{ikt}$ converges
to $F$ in ${\frak G}'_{\{\beta\}}$-topology.
\begin{cor} \
In order that an $m$-harmonic in $K$ function $u(r, t)$ have a boundary value in the space
${\frak G}'_{\{\beta\}}$, it is necessary and sufficiant that the condition (14)
be satisfied.
\end{cor}
\newpage
%\bigskip
%\begin{center}{REFERENCES} \end{center}
%\smallskip

\vspace*{3mm}
Institute of Mathematics \\
National Academy of Sciences of Ukraine \\
3 Tereshchenkivs'ka \\
Kyiv 01601, Ukraine \\
E-mail: imath@horbach.kiev.ua, sergiy.torba@gmail.com
\end{document}